\documentclass[letter,10pt,conference]{IEEEtran}

 \usepackage[bookmarks=false]{hyperref}
\ifCLASSINFOpdf
\else
\fi
\usepackage{relsize}
\usepackage{cite}
\usepackage{amsfonts}
\usepackage{mathrsfs}
\usepackage{bbm}
\usepackage{pifont}
\usepackage{amsfonts}
\usepackage{amsmath, amsthm}
\usepackage{tabularx}
\usepackage{listings}
\usepackage{graphicx}
\usepackage{float}
\usepackage{amssymb}
\usepackage{array}
\usepackage{fancyhdr}
\usepackage{cases}
 \usepackage{supertabular}
\usepackage{url}
\usepackage{fancyhdr}

\usepackage{psfrag}
\usepackage{color}
\usepackage{bbding}

\newcommand{\bP}[1]{{\mathbb{P}}\left[{#1}\right]}
\newcommand{\bE}[1]{{\mathbb{E}}\left[{#1}\right]}

\newcommand{\1}[1]{{\bf 1}\left[#1\right]}       % indicator 1{...}

\newcommand{\fsquare}{\vrule height6pt width7pt depth1pt}   % filled square
\newcommand{\myproof}{{\hfill \\ \bf Proof. \ }}           % Proof
\newcommand{\myendpf}{\hfill\fsquare \\[0.1in]}             % end of proof

\def\centerhack#1{\hbox to 0pt{\hss\footnotesize #1\hss}}
\def\centerhackn#1{\hbox to 0pt{\hss #1\hss}}
\def\dchack#1{\vbox to 0pt{\vss{\hbox to 0pt{\hss#1\hss}}\vss}}

\newtheorem{theorem}{Theorem}[section]

\newtheorem{lemma}[theorem]{Lemma}

\newtheorem{corollary}[theorem]{Corollary}

\begin{document}

\title{On the Eschenauer-Gligor key pre-distribution scheme 
       under on-off communication channels: \\
        The absence of isolated nodes\\
        (Extended version)}

\author{ \IEEEauthorblockN{Armand M. Makowski}
\IEEEauthorblockA{Department of Electrical and Computer Engineering, \\
                                   and the Institute for Systems Research, University of Maryland,\\
                                   College Park, MD 20742 USA \\
                                   E-mail: armand@isr.umd.edu.}
\and \IEEEauthorblockN{Osman Ya\u{g}an}
\IEEEauthorblockA{CyLab and Department
of ECE\\
Carnegie Mellon University \\
Moffett Field, CA 94035 USA \\
Email: oyagan@ece.cmu.edu.}
}

\maketitle \thispagestyle{plain} \pagestyle{plain}

\iffalse

\makeatletter
\let\old@ps@headings\ps@headings
\let\old@ps@IEEEtitlepagestyle\ps@IEEEtitlepagestyle
\def\confheader#1{%
  % for all pages except the first
  \def\ps@headings{%
    \old@ps@headings%
    \def\@oddhead{\strut\hfill#1\hfill\strut}%
    \def\@evenhead{\strut\hfill#1\hfill\strut}%
  }%
  % for the first page
  \def\ps@IEEEtitlepagestyle{%
    \old@ps@headings%
    \def\@oddhead{\strut\hfill#1\hfill\strut}%
    \def\@evenhead{\strut\hfill#1\hfill\strut}%
  }%
  \ps@headings%
} \makeatother

\confheader{%
  \large This paper was submitted to IEEE ISIT in
January 2014 (under review now). }

\markboth{a}{b}

\fancyhead[C]{\large This paper was submitted to IEEE ISIT in
January 2014 (under review now).}

\fi

%\markboth{} {\large This paper has been accepted at IEEE ISIT 2014
%to be held in June--July 2014.}

\begin{abstract}
\boldmath 
We consider the Eschenauer-Gligor key pre-distribution scheme 
under the condition of partial visibility with i.i.d. on-off 
links between pairs of nodes.
This situation is modeled as the intersection
of two random graphs, namely a random key graph
and an Erd\H{o}s-R\'enyi (ER) graph. For this class of composite random graphs we
give various improvements
on a recent result by Ya\u{g}an  \cite{YaganRKGER} concerning
zero-one laws for the absence of isolated nodes.
\end{abstract}

%\vspace{10pt}
%
%\begin{IEEEkeywords}
%key predistribution, node degree, random graph, random intersection
%graph, random key graph, security, topological properties, wireless
%sensor networks.
% \end{IEEEkeywords}


 \section{Introduction}
\label{sec:Introduction}

By now there exists already a large literature discussing various
performance aspects of random key predistribution schemes
in wireless sensor networks (WSNs); see
\cite{CamtepeYener, PerrigStankovicWagner, SunHe, WangAtteburyRamamurthy, XiaoRayiSunDuHuGallowaySurvey}.
However, starting with the scheme of Eschenauer and Gligor
\cite{EschenauerGligor},
much of the work to date has been carried out 
under the {\em full} visibility assumption whereby sensor
nodes are all within communication range of each other. 
While the full visibility assumption is certainly at
odds with the wireless nature of the communication medium supporting WSNs, 
this simplification makes it possible to focus solely 
on how the randomization mechanism affects performance 
in the best of circumstances,
i.e., when wireless communication is not a bottleneck. 
A common criticism of this line of work is that by
disregarding the unreliability of the wireless links, 
the resulting dimensioning guidelines are likely to be overly
optimistic, if not irrelevant. 
In practice, nodes will have fewer neighbors since some of the communication
links may be impaired.

In a recent paper \cite{YaganRKGER}, Ya\u{g}an studied
the Eschenauer-Gligor key pre-distribution scheme 
under the condition of {\em partial} visibility with i.i.d. on-off 
links between pairs of nodes.
This situation was modeled as the {\em intersection}
of two random graphs, namely a random key graph
\cite{BlackburnGerke, DiPietroManciniMeiPanconesiRadhakrishnan2008, YaganThesis,
YaganMakowskiISIT2008, YaganMakowskiISIT2009, YaganMakowskiConnectivity}
and an Erd\H{o}s-R\'enyi (ER) graph \cite{Bollobas, JansonLuczakRucinski}:
With $n$ nodes in the network, 
the Eschenauer-Gligor scheme 
with key rings of size $K$ drawn from a pool of $P$ distinct keys
($K < P$) gives rise to the random key graph $\mathbb{K}(n;\theta)$
(where we have set $\theta = (K,P)$) -- Let $q(\theta)$ denote
the probability (\ref{eq:q_theta}) that a link does not exist between two nodes
in $\mathbb{K}(n;\theta)$.
The communication model between nodes
corresponds to an Erd\H{o}s-R\'enyi (ER) graph $\mathbb{G}(n;\alpha)$
with link probability $\alpha$ (in $[0,1]$).
Under a natural independence assumption,
the graph of interest is the graph 
$\mathbb{K\cap G}(n;\theta,\alpha)$ whose edge set
is the intersection of the edge sets of the random graphs
$\mathbb{K}(n;\theta)$ and $\mathbb{G}(n;\alpha)$.
See Section \ref{sec:Model} for more details concerning the model and the notation in use.
            
In \cite{YaganRKGER}
the following zero-one law for the absence of isolated nodes
was established: If the parameters are scaled 
with the number $n$ of nodes in such a way that
\begin{equation}
\alpha_n 
\left (1 - q (\theta_n) \right ) \sim c \frac{\log n}{n}
\label{eq:YaganResultCondition}
\end{equation}
for some $c> 0$, then it holds that
\begin{eqnarray}
\lefteqn{
\lim_{n \rightarrow \infty}
\bP{
\begin{array}{c}
\mathbb{K} \cap \mathbb{G} (n; \theta_n, \alpha_n)
\mbox{~contains }  \\
\mbox{~no isolated nodes}  \\
\end{array}
}
} & &
\nonumber \\
&=& \left \{
\begin{array}{ll}
0 & \mbox{if ~$0<c<1$} \\
  &     \\
1 & \mbox{if ~$1<c$} \\
\end{array}
\right .
\label{eq:YaganResult}
\end{eqnarray}
provided the limit $\lim_{n \rightarrow \infty} \alpha_n \log n $
{\em exists} in $[0,\infty]$.
%This improved on the work of Yi et al.
%\cite{YiWanLinHuang} who established
%the zero-law under the restrictive assumption that
%$\lim_{n \rightarrow \infty} \alpha_n \log n = \infty$.

In this short paper, we improve on this result 
in two different directions which are now briefly described. 
Precise statements are available in Section \ref{sec:MainResults}:

(i) We show that 
            the existence of a limit for the sequence
            $\{ \alpha_n \log n, \ n=2,3, \ldots \}$ 
            is not needed to ensure the zero-one law (\ref{eq:YaganResult}) 
            under (\ref{eq:YaganResultCondition}).
            In fact, this result was already contained
            in the earlier result of Ya\u{g}an \cite{YaganRKGER}, and is an
            easy consequence of the Principle of Subsubsequences
             \cite{JansonLuczakRucinski}.

(ii) We partially strengthen the result of Ya\u{g}an \cite{YaganRKGER} 
               by establishing a zero-one law when the scaling is done according to
\begin{equation}
\alpha_n ( 1 - q(\theta_n) ) = \frac{\log n + \gamma_n }{n},
\quad n=1,2, \ldots
\label{eq:MakowskiYaganResultCondition}
\end{equation}
            for some deviation function 
            $\gamma : \mathbb{N}_0 \rightarrow \mathbb{R}$.
            This is done under mild conditions on the scaling $\{ \alpha_n, \ n=1,2, \ldots \}$.
            The class of scalings satisfying (\ref{eq:YaganResultCondition}) 
            is easily seen to be contained in the class of scalings
            governed by (\ref{eq:MakowskiYaganResultCondition}).
            The proof uses the method of first and second moments applied 
            to the number of isolated nodes --
            This approach is presented in Section \ref{sec:MethodFirstSecond} where expressions for the needed moments
            are given; see Appendix \ref{sec:Appendix} for detailed calculations.
            The asymptotics of the first moment
            are derived in Section \ref{sec:BehaviorFirstMoment} in terms of a \lq\lq zero-infinity" law.
            The bounds for applying the method of second moment are derived in Section \ref{sec:Bounds}.
            The proof of the zero-law under the scaling (\ref{eq:MakowskiYaganResultCondition}) is completed in
            Section \ref{sec:ProofPartialVeryStrongZeroLaw+NodeIsolation},
            Section \ref{AlongSubsequences} and Section \ref{sec:CompleteProof}.
            
         The material of this paper appeared in the 
         Proceedings of the 53rd Annual Allerton Conference on Communication, Control, and Computing, Monticello (IL)
         \cite{MakowskiYagan_Allerton}.    
            
\section{The model}
\label{sec:Model}

All limiting statements, including asymptotic equivalences, are understood with
the number $n$ of sensor nodes going to infinity. The random
variables (rvs) under consideration are all defined on the same
probability triple $(\Omega, {\cal F}, \mathbb{P})$. Probabilistic
statements are made with respect to this probability measure
$\mathbb{P}$, and we denote the corresponding expectation operator
by $\mathbb{E}$. The indicator function of an event $E$ is
denoted by $\1{E}$. 
For any discrete set $S$ we write $|S|$ for its cardinality. 

\subsection{The Eschenauer-Gligor scheme}
\label{subsec:EG}

The Eschenauer-Gligor scheme is characterized by three parameters,
which are held fixed throughout this section, namely the number $n$ of nodes, 
the size $P$ of the key pool and the size $K$ of each key ring 
with $K < P$. 
To lighten the notation we often group the integers $P$
and $K$ into the ordered pair $\theta \equiv (K,P)$.

Nodes are labelled $i=1, \ldots , n$.
For each $i=1, \ldots , n$, let $K_i (\theta)$ denote the
random set of $K$ distinct keys assigned to node $i$ before network
deployment. 
According to the Eschenauer-Gligor scheme,
if after deployment,
two nodes, say $i$ and $j$, are within communication range of
each other, they can establish a secure link provided 
their key rings have at least one key in common.

We can think
of $K_i(\theta)$ as an $\mathcal{P}_{K} $-valued rv where
$\mathcal{P}_{K} $ denotes the collection of all subsets of 
$\{ 1, \ldots , P \}$ which contain exactly $K$ elements -- Obviously, 
we have $|\mathcal{P}_{K} | =  {P \choose K}$. The rvs $K_1(\theta),
\ldots , K_n(\theta)$ are assumed to be {\em i.i.d.} rvs, each of
which is {\em uniformly} distributed over $\mathcal{P}_{K}$ with
\begin{equation}
\bP{ K_i(\theta) = S } 
= {P \choose K} ^{-1}, 
\quad 
\begin{array}{c}
i=1, \ldots , n \\
S \in \mathcal{P}_{K} .\\
\end{array}
\label{eq:EG_KeyDistrbution1}
\end{equation}
This corresponds to selecting keys
randomly and {\em without} replacement from the key pool.

For future reference, for any subset $R$ of $\{ 1, \ldots , P \}$
we find it convenient to write
\begin{equation}
v(\theta;R)
=
\left \{
\begin{array}{ll}
\frac{{P- |R| \choose K}}{{P \choose K}} &
\mbox{if~ $|R| \leq P-K$} \\
  &  \\
0 &
\mbox{if~ $P-K < |R|$.} \\
\end{array}
\right .
\label{eq:V}
\end{equation}
Since
$v(\theta;R)$ depends on $R$ only through its cardinality $|R|$, 
sometimes we shall also write $v(\theta;|R|)$ in place of $v(\theta;R)$. 
It is a simple matter to check that
\begin{equation}
\bP{ K_i(\theta) \cap R = \emptyset } = v(\theta;R),
\quad i=1, \ldots , n.
\label{eq:Probab_Keyring_Does_Not_Intersect_R}
\end{equation}

\subsection{Random key graphs}

Under full visibility, the Eschenauer-Gligor scheme gives rise to a random graph which
we now describe:
Distinct nodes $i$ and $j$ are said to be $K$-adjacent,
written $i \sim_{K} j$, if 
their key rings have at least one key in common.
Thus,
\begin{equation}
i \sim_K j \quad \mbox{iff} \quad 
K_i (\theta) \cap K_j (\theta) \neq \emptyset ,
\label{eq:Adjacency_RKG}
\end{equation} 
and an undirected link is assigned between nodes $i$ and $j$. 
This notion of adjacency defines the {\em random key graph} 
$\mathbb{K}(n; \theta )$ on the vertex set $\{ 1, \ldots , n\}$.

For distinct $i,j =1, \ldots , n$, it is
a simple matter to check from (\ref{eq:Probab_Keyring_Does_Not_Intersect_R}) that
\begin{equation}
\bP{ K_i (\theta) \cap K_j (\theta) = \emptyset } = q (\theta)
\end{equation}
with
\begin{equation}
q (\theta) = \left \{
\begin{array}{ll}
0 & \mbox{if~ $P <2K$} \\
  &                 \\
\frac{{P-K \choose K}}{{P \choose K}} & \mbox{if~ $2K \leq P$.}
\end{array}
\right . 
\label{eq:q_theta}
\end{equation}
Note that $q(\theta) = v(\theta, K)$.
It is plain that
\begin{equation}
\bP{i \sim_{K} j} = 1 -q(\theta) 
\label{eq:Edge_Prob_RKG}
\end{equation}
so that the probability of edge occurrence between any two nodes is
equal to $1-q(\theta)$. 

\subsection{ER graphs as a simple communication model}

To account for the
possibility that communication links between nodes may not be
available, we assume a simple communication model
that consists of independent communication channels,
each of which can be either on or off. 
Thus, with $\alpha$ in $[0,1]$, 
let $\{B_{ij}(\alpha), \ 1 \leq i < j \leq n\}$ 
denote i.i.d. $\{0, 1\}$-valued rvs with
success probability $\alpha$. For convenience we also introduce
the $\{0, 1\}$-valued rvs
$\{B_{ji}(\alpha), \ 1 \leq i < j \leq n\}$
by setting
\[
B_{ji} (\alpha) = B_{ij} (\alpha),
\quad 1 \leq i < j \leq n.
\]

The channel between nodes $i$ and $j$ is
available (equivalently, up) if $B_{ij}(\alpha) =1$ 
with probability $\alpha$, and unavailable (equivalently, down) if $B_{ij}(\alpha) = 0$ 
with complementary probability $1-\alpha$.
Distinct nodes $i$ and $j$ are said to be B-adjacent, 
written $i \sim_{B} j$, if $B_{ij}(\alpha) = 1$.
The notion of B-adjacency defines
the standard ER graph $\mathbb{G}(n;\alpha)$ on the vertex set 
$\{ 1, \ldots , n \}$. Obviously,
\[
\bP{ i \sim_{B} j} = \alpha.
\]

\subsection{Intersecting the graphs}

The random graph model studied here is obtained 
by {\em intersecting} 
the random key graph
$\mathbb{K}(n;\theta)$ with the
ER graph $\mathbb{G}(n;\alpha)$: The distinct nodes $i$
and $j$ are now said to be adjacent, written $i \sim j$, 
if and only if they are both K-adjacent and B-adjacent, namely
\begin{equation}
i \sim j \quad \mbox{iff} \quad
\begin{array}{c}
K_i (\theta) \cap K_j (\theta) \neq \emptyset \\
\mbox{and} \\
B_{ij}(\alpha)=1.\\
\end{array}
\label{eq:Adjacency_Intersection}
\end{equation}
The resulting undirected random graph defined on the vertex
set $\{1, \ldots, n\}$ through this notion of adjacency is denoted
$\mathbb{K\cap G}(n;\theta,\alpha)$.

Throughout, the collections of rvs 
$\{ K_1 (\theta), \ldots , K_n (\theta) \}$ 
and $\{B_{ij}(\alpha), \ 1 \leq i < j \leq n\}$ 
are assumed to be {\em independent}, 
in which case the probability of edge occurrence 
in $\mathbb{K\cap G}(n;\theta,\alpha)$ is given by
\begin{eqnarray}
\bP{i \sim j} 
= \bP{ i \sim_{K} j } \bP{ i \sim_{B} j }
= p(\theta,\alpha)
\label{eq:Edge_Prob_Intersection}
\end{eqnarray}
where we have set
\begin{equation}
p(\theta,\alpha) 
= \alpha ( 1 - q (\theta) ).
\label{eq:p(theta,alpha)}
\end{equation}
Finally, to simplify the notation, we set
\[
P_n (\theta, \alpha)
=
\bP{
\begin{array}{c}
\mathbb{K} \cap \mathbb{G} (n; \theta_n, \alpha_n)
\mbox{~contains }  \\
\mbox{~no isolated nodes}  \\
\end{array}
}.
\]

\section{The main results}
\label{sec:MainResults}

To fix the terminology, we refer to any pair of mappings
$K,P: \mathbb{N}_0 \rightarrow \mathbb{N}_0$ as a scaling
(for random key graphs)
provided the natural conditions
\begin{equation}
K_n < P_n, \quad n=1,2, \ldots . \label{eq:ScalingDefn}
\end{equation}
are satisfied. 
Similarly, any mapping $\alpha: \mathbb{N}_0 \rightarrow [0,1]$ defines
a scaling for ER graphs.

The terminology of strong and very strong zero-one laws parallels
the one introduced in the survey papers
\cite[Section IV, p. 1070]{HanMakowski_JSAC}
\cite{MakowskiHan_NOW}.
The first result gives a very strong
one-law for the absence of isolated nodes under minimal assumptions;
its proof is given in Section \ref{sec:BehaviorFirstMoment}.

\begin{theorem}
{\sl
Consider scalings $K,P: \mathbb{N}_0 \rightarrow \mathbb{N}_0$
and $\alpha: \mathbb{N}_0 \rightarrow [0,1]$ such that
\begin{equation}
\alpha_n ( 1 - q(\theta_n) ) = \frac{\log n + \gamma_n }{n},
\quad n=1,2, \ldots
\label{eq:ScalingCondition+VeryStrong}
\end{equation}
for some deviation function $\gamma : \mathbb{N}_0 \rightarrow \mathbb{R}$.
The very strong one-law
\begin{eqnarray}
\lim_{n \rightarrow \infty} P_n (\theta_n, \alpha_n)
= 1
\label{eq:VeryStrongOneLaw+NodeIsolation}
\end{eqnarray}
holds whenever
\begin{equation}
\lim_{n\rightarrow \infty} \gamma_n = \infty .
\label{eq:ConditionForVeryStrongOneLaw+NodeIsolation}
\end{equation}
}
\label{thm:VeryStrongOneLaw+NodeIsolation}
\end{theorem}

It is noteworthy that Theorem \ref{thm:VeryStrongOneLaw+NodeIsolation} applies
to the constant parameter case, yielding a result similar to the one available for many
classes of random graphs, e.g., ER graphs \cite{Bollobas, JansonLuczakRucinski},
geometric random graphs \cite{GuptaKumar} and
random key graphs \cite{YaganMakowskiCISS2010}.
The proof is straightforward and is omitted in the interest of brevity.

\begin{corollary}
{\sl
With $\alpha$ in $(0,1]$ and positive integers $K$ and $P$ such that $K < P$, we always have
$\lim_{n \rightarrow \infty} P_n (\theta, \alpha) = 1$ provided $\alpha (1-q(\theta) ) > 0$.
}
\end{corollary}

\myproof
We can write
\[
\alpha (1-q(\theta) ) = \frac{\log n + \gamma_n }{n},
\quad n=1,2, \ldots
\]
with deviation function $\gamma: \mathbb{N}_0 \rightarrow \mathbb{R}$ given by
\[
\gamma_n = n \alpha (1-q(\theta) ) - \log n,
\quad n=1,2, \ldots
\]
The desired conclusion is a simple consequence of Theorem \ref{thm:VeryStrongOneLaw+NodeIsolation} 
as we note that $\lim_{n \rightarrow \infty} \gamma_n = \infty$ under the condition
$\alpha (1-q(\theta) ) > 0$.
\myendpf

While no additional condition are needed in
Theorem \ref{thm:VeryStrongOneLaw+NodeIsolation},
the corresponding zero-law does require growth conditions
on the scaling $\alpha: \mathbb{N}_0 \rightarrow [0,1]$.
\begin{theorem}
{\sl
Consider scalings $K,P: \mathbb{N}_0 \rightarrow \mathbb{N}_0$
and $\alpha: \mathbb{N}_0 \rightarrow [0,1]$ such that
(\ref{eq:ScalingCondition+VeryStrong}) holds
for some deviation function $\gamma : \mathbb{N}_0 \rightarrow \mathbb{R}$.
The very strong zero-law
\begin{eqnarray}
\lim_{n \rightarrow \infty} P_n (\theta_n, \alpha_n) = 0
\label{eq:VeryStrongZeroLaw+NodeIsolation}
\end{eqnarray}
holds whenever
\begin{equation}
\lim_{n\rightarrow \infty} \gamma_n = -\infty 
\label{eq:ConditionForVeryStrongZeroLaw+NodeIsolation}
\end{equation}
provided either 
\begin{equation}
\limsup_{n\rightarrow \infty} \alpha_n \log n < \infty,
\label{eq:AdditionalConditionForVeryStrongZeroLaw+NodeIsolation}
\end{equation}
or
\begin{equation}
\limsup_{n\rightarrow \infty} \alpha_n \log n = \infty
\quad \mbox{with} \quad \limsup_{n\rightarrow \infty} \alpha_n < 1.
\label{eq:AdditionalConditionForVeryStrongZeroLaw+NodeIsolation2}
\end{equation}
}
\label{thm:PartialVeryStrongZeroLaw+NodeIsolation}
\end{theorem}

A proof of Theorem \ref{thm:PartialVeryStrongZeroLaw+NodeIsolation} is developed 
through Sections \ref{sec:MethodFirstSecond} to \ref{sec:ProofPartialVeryStrongZeroLaw+NodeIsolation}.
The additional growth conditions
(\ref{eq:AdditionalConditionForVeryStrongZeroLaw+NodeIsolation})-(\ref{eq:AdditionalConditionForVeryStrongZeroLaw+NodeIsolation2})
can be dropped when restricting attention to the scalings used by
Ya\u{g}an \cite{YaganRKGER}.

\begin{theorem}
{\sl 
Consider scalings $K,P: \mathbb{N}_0 \rightarrow \mathbb{N}_0$
and $\alpha: \mathbb{N}_0 \rightarrow [0,1]$ such that
\begin{equation}
\alpha_n ( 1 - q(\theta_n) ) \sim c ~ \frac{\log n }{n}
\label{eq:ScalingCondition+Strong}
\end{equation}
for some $c>0$. 
Then, the strong zero-one law
\begin{eqnarray}
\lim_{n \rightarrow \infty} P_n (\theta_n, \alpha_n)
=
\left \{
\begin{array}{ll}
0 & \mbox{if ~$0<c<1$} \\
  &     \\
1 & \mbox{if ~$1<c$} \\
\end{array}
\right .
\label{eq:StrongOneLaw+NodeIsolation}
\end{eqnarray}
holds.
} 
\label{thm:StrongOneLaw+NodeIsolation}
\end{theorem}

\myproof
Consider scalings $K,P: \mathbb{N}_0 \rightarrow \mathbb{N}_0$
and $\alpha: \mathbb{N}_0 \rightarrow [0,1]$ such that
(\ref{eq:ScalingCondition+Strong})
holds for some $c>0$. 
This can be rewritten in equivalent form as
\begin{equation}
\alpha_n ( 1 - q(\theta_n) ) = c_n \cdot \frac{\log n }{n},
\quad n=1,2, \ldots
\label{eq:ScalingCondition+Strong+Equivalent}
\end{equation}
where the sequence $c: \mathbb{N}_0 \rightarrow \mathbb{R}_+$
satisfies $\lim_{n \rightarrow \infty} c_n = c$.
It is then plain that
(\ref{eq:ScalingCondition+VeryStrong}) automatically holds with deviation function 
$\gamma : \mathbb{N}_0 \rightarrow \mathbb{R}$ given by
\[
\gamma_n = (c_n -1 ) \log n,
\quad n=1,2, \ldots
\]

When $c > 1$, we have $\lim_{n \rightarrow \infty} \gamma_n = \infty$ and
Theorem \ref{thm:VeryStrongOneLaw+NodeIsolation} gives the one-law 
(\ref{eq:VeryStrongOneLaw+NodeIsolation}), hence the one-law
part of (\ref{eq:StrongOneLaw+NodeIsolation}) holds.
On the other hand, with $0 < c < 1$, 
$\lim_{n \rightarrow \infty} \gamma_n = - \infty$ and
Theorem \ref{thm:PartialVeryStrongZeroLaw+NodeIsolation}
yields  the zero-law (\ref{eq:VeryStrongZeroLaw+NodeIsolation}),
hence the zero-law part of (\ref{eq:StrongOneLaw+NodeIsolation}),
if the additional conditions
(\ref{eq:AdditionalConditionForVeryStrongZeroLaw+NodeIsolation}) or
(\ref{eq:AdditionalConditionForVeryStrongZeroLaw+NodeIsolation2})
hold.
We now show that this additional condition is superfluous for the zero-law
to hold; this is a consequence of the Principle of Subsubsequences
\cite{JansonLuczakRucinski} -- In what follows a subsequence $k \rightarrow n_k$
is simply any non-decreasing mapping
$\mathbb{N}_0 \rightarrow \mathbb{N}_0: k \rightarrow n_k$ such that
$\lim_{k \rightarrow \infty} n_k = \infty$:

A careful inspection of the arguments given by Ya\u{g}an
\cite[Thm. 3.1, p. 3824]{YaganRKGER} shows that the 
result also holds along subsequences:
Specifically, consider scalings $K,P: \mathbb{N}_0 \rightarrow \mathbb{N}_0$
and $\alpha: \mathbb{N}_0 \rightarrow [0,1]$ such that
(\ref{eq:ScalingCondition+Strong})
holds for some $c$ in $(0,1)$.
Then, for any subsequence $k \rightarrow n_k$, we have
\begin{eqnarray}
\lim_{k \rightarrow \infty}
P_{n_k} (\theta_{n_k},\alpha_{n_k} )
= 0
\label{eq:StrongZeroLaw+NodeIsolationAlongSubsequence}
\end{eqnarray}
whenever the limit $\lim_{k \rightarrow \infty} \alpha_{n_k} \log n_k$ 
exists in $[0,\infty]$.

The sequence
\begin{equation}
\{ P_n (\theta_n,\alpha_n ) , \ n=2,3, \ldots \}
\label{eq:TheSequence}
\end{equation}
is a bounded sequence with {\em all} its accumulation points in $[0,1]$.
Let $P$ be {\em any} accumulation point of the sequence.
By definition, there exists a subsequence $k \rightarrow n_k$ such that
\begin{equation}
\lim_{k \rightarrow \infty}
P_{n_k} (\theta_{n_k},\alpha_{n_k} ) = P.
\label{eq:AlongSubsequence1}
\end{equation}
Although the sequence
$\{ \alpha_{n_k} \log n_k, \ k=1,2, \ldots \}$ may not converge,
there must exist a further subsequence
$\ell \rightarrow k_\ell$ such that the limit
$ \lim_{\ell \rightarrow \infty} \alpha_{n_{k_\ell}} \log n_{k_\ell}$
does exist in $[0,\infty]$. 

Taking 
(\ref{eq:AlongSubsequence1}) along that subsequence we find
\[
\lim_{\ell \rightarrow \infty} 
P_{n_{k_\ell}} (\theta_{n_{k_\ell}},\alpha_{n_{k_\ell}} ) = P,
\]
whence $P=0$ by virtue of (\ref{eq:StrongZeroLaw+NodeIsolationAlongSubsequence}).
The {\em bounded} sequence (\ref{eq:TheSequence})
thus admits $P=0$ as its {\em unique} accumulation point,
and is therefore convergent with limit
\[
\lim_{n \rightarrow \infty}
P_n (\theta_n,\alpha_n )
= 0
\]
regardless of whether the
sequence $\{ \alpha_n \log n , \ n=1,2, \ldots \}$ 
has a limit in $[0,\infty]$.
\myendpf

\section{The method of first and second moments}
\label{sec:MethodFirstSecond}

Theorem \ref{thm:VeryStrongOneLaw+NodeIsolation} 
and
Theorem \ref{thm:PartialVeryStrongZeroLaw+NodeIsolation} 
will be established by 
the method of first and second moments \cite[p.  55]{JansonLuczakRucinski}
applied to the number of isolated nodes. 
Fix $n=2,3, \ldots $ and
consider positive integers $K$ and $P$ such that $K < P$, 
and scalar $\alpha$ in $[0,1]$.

\subsection{Counting isolated nodes}

The number of isolated nodes in
$\mathbb{K \cap G}(n;\theta,\alpha)$ is given by
\[
I_n (\theta,\alpha) = \sum_{i=1}^n \chi_{n,i}(\theta,\alpha)
\]
where  for each $i=1,2, \ldots , n$, we write
\[
\chi_{n,i}(\theta,\alpha) 
= 
\1{ \mbox{Node $i$ is isolated in $\mathbb{K \cap G}(n;\theta, \alpha)$ } }.
\]
It is a simple matter to check that
\begin{equation}
\chi_{n,i}(\theta,\alpha)
= \prod_{j=1, \ j \neq i}^n 
\left (
1 - B_{ij}(\alpha) \eta_{ij}(\theta)
\right )
\label{eq:Defn+CHI}
\end{equation}
with indicator rvs
\begin{equation}
\eta_{ij}(\theta) 
= \1{ K_i (\theta) \cap K_j(\theta) \neq \emptyset },
\quad 
\begin{array}{c}
i \neq j \\
i,j=1, \ldots , n. \\
\end{array}
\end{equation}
The random graph $\mathbb{K \cap G}(n;\theta,\alpha)$ has no isolated
nodes if and only if $I_n (\theta,\alpha) = 0$, and 
the key relation
\[
P_n(\theta,\alpha)
= \bP{ I_n(\theta,\alpha) = 0 }
\]
follows.

This equivalence is exploited with the help of two standard bounds
based on first and second moments:
The easy bound 
\begin{equation}
1 - \bE{ I_n (\theta,\alpha) } \leq \bP{  I_n (\theta,\alpha) = 0 }
\label{eq:BoundViaFirstMoment}
\end{equation}
gives rise to the method of first moment
\cite[Eqn. (3.10), p. 55]{JansonLuczakRucinski},
while the method of second moment 
\cite[Remark 3.1, p. 55]{JansonLuczakRucinski} 
has its starting point in the inequality
\begin{equation}
\bP{  I_n (\theta,\alpha) = 0 } 
\leq 1 - \frac{ \left ( \bE{ I_n (\theta,\alpha)} \right )^2}
                     { \bE{ I_n (\theta,\alpha) ^2} }. 
\label{eq:BoundViaSecondMoment}
\end{equation}

\subsection{Evaluating moments}

The rvs $\chi_{n,1}(\theta,\alpha), \ldots , \chi_{n,n} (\theta,\alpha)$ 
being exchangeable, we readily get
\begin{equation}
\bE{ I_n (\theta,\alpha)} = n \bE{ \chi_{n,1} (\theta,\alpha) }
\label{eq:FirstMomentExpression}
\end{equation}
and
\begin{eqnarray}
\bE{ I_n (\theta,\alpha)^2 } 
&=&
 n \bE{ \chi_{n,1} (\theta,\alpha) }
\nonumber \\
& & +  n(n-1) \bE{ \chi_{n,1} (\theta,\alpha) \chi_{n,2} (\theta,\alpha) } .
\label{eq:SecondMomentExpression}
\nonumber
\end{eqnarray}
This last expression is an easy consequence
of the binary nature of the rvs involved. It then follows that
\begin{eqnarray}
\frac{ \bE{ I_n (\theta,\alpha)^2 }}{ \left (  \bE{ I_n (\theta,\alpha) } \right )^2 } 
&=& 
\frac{ 1}{ \bE{ I_n (\theta,\alpha) } }
\label{eq:SecondMomentRatio} \\
& & + \frac{n-1}{n} \cdot \frac{\bE{ \chi_{n,1} (\theta,\alpha)
\chi_{n,2} (\theta,\alpha) }}
     {\left (  \bE{ \chi_{n,1} (\theta,\alpha) } \right )^2 }.
\nonumber 
\end{eqnarray}

With (\ref{eq:Defn+CHI}) as point of departure, 
expressions are easily obtained for the needed moments
$\bE{ \chi_{n,1} (\theta,\alpha) }$
and
$\bE{ \chi_{n,1} (\theta,\alpha) \chi_{n,2} (\theta,\alpha) }$;
calculations are given in Appendix \ref{sec:Appendix}
for sake of completeness:
In the notation (\ref{eq:p(theta,alpha)}), we have
\begin{equation}
\bE{ \chi_{n,1} (\theta,\alpha) }
=
\left ( 1 - p(\theta,\alpha) \right )^{n-1},
\label{eq:FirstMomentCHI}
\end{equation}
whence
\begin{equation}
\bE{ I_n(\theta,\alpha ) }
= n \left ( 1 - p(\theta,\alpha) \right )^{n-1} .
\label{eq:FirstMoment2}
\end{equation}
We also show that  
\begin{eqnarray}
\lefteqn{
\bE{ \chi_{n,1} (\theta,\alpha) \chi_{n,2} (\theta,\alpha) }
} & &
\nonumber \\
&=&
\bE{ (1- \alpha \eta_{12}(\theta) ) Z (\theta,\alpha)^{n-2} }
\label{eq:CrossMoment}
\end{eqnarray}
where the auxiliary rv $Z(\theta,\alpha)$ is given by
\begin{equation}
Z (\theta,\alpha)
=
\left ( 1 - p(\theta,\alpha ) \right )^2
\left ( 1 + \widetilde Z (\theta,\alpha)  \right )
\label{eq:Z}
\end{equation}
with
\begin{eqnarray}
\lefteqn{
\widetilde Z (\theta,\alpha)
} & &
\label{eq:TildeZ} \\
&=&
\frac{ \alpha^2 }{ \left ( 1 - p(\theta,\alpha ) \right )^2 }
\cdot 
\left ( v\left(\theta; K_1(\theta) \cup K_2(\theta) \right) - q(\theta)^2 \right ).
\nonumber
\end{eqnarray} 

\section{Behavior of the first moment}
\label{sec:BehaviorFirstMoment}

The proof of Theorem \ref{thm:VeryStrongOneLaw+NodeIsolation} passes through a characterization of the 
behavior of the first moment given in the 
following \lq\lq zero-infinity" law -- Note its \lq\lq analogy" 
with Theorem \ref{thm:VeryStrongOneLaw+NodeIsolation}.

\begin{lemma}
{\sl
Consider scalings $K,P: \mathbb{N}_0 \rightarrow \mathbb{N}_0$
and $\alpha: \mathbb{N}_0 \rightarrow [0,1]$ such that
(\ref{eq:ScalingCondition+VeryStrong}) holds
for some deviation function $\gamma : \mathbb{N}_0 \rightarrow \mathbb{R}$.
It is always the case that
\begin{eqnarray}
\lim_{n \rightarrow \infty} \bE{ I_n (\theta_n,\alpha_n) }
= \left \{
\begin{array}{ll}
\infty & \mbox{if ~$\lim_{n\rightarrow \infty} \gamma_n = -\infty$} \\
  &     \\
0 & \mbox{if ~$\lim_{n\rightarrow \infty} \gamma_n = \infty$.} \\
\end{array}
\right .
\label{eq:VeryStrongZeroOneLaw+FirstMoment}
\end{eqnarray}
}
\label{lem:VeryStrongZeroOneLaw+FirstMoment}
\end{lemma}
Before establishing this result, we note that the proof of
Theorem \ref{thm:VeryStrongOneLaw+NodeIsolation}
is now straightforward:
The bound (\ref{eq:BoundViaFirstMoment}) yields
$\lim_{n\to \infty} \bP{I_n(\theta_n,\alpha_n) = 0} = 1$ 
whenever $\lim_{n \to \infty} \bE{ I_n (\theta_n,\alpha_n) }= 0$,
as this is the case under the condition
(\ref{eq:ConditionForVeryStrongOneLaw+NodeIsolation}) by virtue of
Lemma \ref{lem:VeryStrongZeroOneLaw+FirstMoment}.
\myendpf

Although the proof of Lemma \ref{lem:VeryStrongZeroOneLaw+FirstMoment} 
is fairly standard, we give some of the details as we need 
to develop some facts that will be used later:
We start with the observation that for $0 \leq x < 1$,
\[
\log (1-x)
= - x - \Psi(x)
\quad \mbox{with} \quad 
\Psi(x) = \int_0^x \frac{t}{1-t} dt.
\]
It is also easy to check that
\begin{equation}
\lim_{x \downarrow 0} \frac{ \Psi(x) }{x^2 } = \frac{1}{2}.
\label{eq:OrderOfPsiAtZero}
\end{equation}

Now consider scalings $K,P: \mathbb{N}_0 \rightarrow \mathbb{N}_0$
and $\alpha: \mathbb{N}_0 \rightarrow [0,1]$ such that
(\ref{eq:ScalingCondition+VeryStrong}) holds
for some deviation function $\gamma : \mathbb{N}_0 \rightarrow \mathbb{R}$.
For each $n=1,2, \ldots$, substitution 
of (\ref{eq:ScalingCondition+VeryStrong})  into
(\ref{eq:FirstMoment2}) yields
\begin{eqnarray}
 \bE{ I_n (\theta_n,\alpha_n) }
 &=&
 n   \left ( 1 - p(\theta_n, \alpha_n)  \right )^{n-1}
 \nonumber \\
 &=&
 n
e^{ (n-1) \log ( 1 - p(\theta_n, \alpha_n)  ) }
 \nonumber \\
 &=&
 n
 e^{-(n-1)
 \left (  p(\theta_n, \alpha_n)  + \Psi( p(\theta_n, \alpha_n) )
 \right )}
 \nonumber \\
 &=&
 n
 e^{- (n-1) \frac{ \log n + \gamma_n }{n} - (n-1) \Psi( p(\theta_n, \alpha_n) ) }
 \nonumber \\
 &=&
n^{ \frac{1}{n}} 
e^{-\frac{n-1}{n} \gamma_n } e^{- (n-1) \Psi( p(\theta_n, \alpha_n) ) }
\label{eq:AA}
\end{eqnarray}
as well as the bound
\begin{equation}
 \bE{ I_n (\theta_n,\alpha_n) } 
\leq
n^{ \frac{1}{n}} 
e^{-\frac{n-1}{n} \gamma_n }.
\label{eq:AB}
\end{equation}

If $\lim_{n \rightarrow \infty} \gamma_n = \infty$, then
$\lim_{n \to \infty} \bE{ I_n (\theta_n,\alpha_n) }= 0$ by virtue
of the inequality (\ref{eq:AB}).
On the other hand, the condition $\lim_{n \rightarrow \infty} \gamma_n 
= -\infty$
already implies
$\lim_{n \rightarrow \infty} 
n^{ \frac{1}{n}}  e^{-\frac{n-1}{n} \gamma_n } = \infty$.
In view of (\ref{eq:AA}),
the desired conclusion $\lim_{n \to \infty} \bE{ I_n (\theta_n,\alpha_n) } = \infty$
then holds if we show
\begin{equation}
\lim_{n \rightarrow \infty}  (n-1) \Psi( p(\theta_n, \alpha_n) )  = 0.
\label{eq:AC}
\end{equation}
To do so, note that 
the condition $\lim_{n \rightarrow \infty} \gamma_n = -\infty$ also
implies $\gamma_n < 0$ for all $n$ sufficiently large, in which case
$\gamma_n = -|\gamma_n|$. On that range
the condition (\ref{eq:ScalingCondition+VeryStrong}) becomes
\[
0 \leq p(\theta_n,\alpha_n) = \frac{ \log n - |\gamma_n|}{n} ,
\]
whence
\begin{equation}
|\gamma_n| \leq \log n
\quad \mbox{and} \quad
p(\theta_n,\alpha_n) \leq \frac{\log n }{n}.
\label{eq:Inequalities}
\end{equation}
Therefore, we must have
\begin{equation}
\lim_{n\rightarrow \infty} p(\theta_n,\alpha_n) = 0
\label{eq:Limit=Zero}
\end{equation}
as well as 
\begin{equation}
\lim_{n\rightarrow \infty} (n-1) p(\theta_n,\alpha_n)^2 = 0.
\label{eq:Limit=Zero2}
\end{equation}
The conclusion (\ref{eq:AC}) is an easy consequence of
these two facts (combined with (\ref{eq:OrderOfPsiAtZero}))
once we note that
\[
(n-1) \Psi( p(\theta_n, \alpha_n) )  
=
(n-1) p(\theta_n, \alpha_n) ^2 
\cdot
\frac{ \Psi( p(\theta_n, \alpha_n) )   }{ p(\theta_n, \alpha_n)^2 }
\]
for all $n=1,2, \ldots$.
\myendpf

\section{Bounds}
\label{sec:Bounds}

The proof of the zero-law relies
on various bounds which we now develop.
Fix $n=2,3, \ldots $ and
consider positive integers $K$ and $P$ such that $K < P$, 
and scalar $\alpha$ in $[0,1]$.

By uninteresting calculations it follows from
(\ref{eq:FirstMomentCHI}), (\ref{eq:CrossMoment}), (\ref{eq:Z}) and (\ref{eq:TildeZ})
that
\begin{eqnarray}
\frac{
\bE{ \chi_{n,1} (\theta,\alpha) \chi_{n,2} (\theta,\alpha) }
}{
\left ( \bE{ \chi_{n,1} (\theta,\alpha) } \right )^2
}
=
\left ( 1 - p(\theta,\alpha) \right )^{-2}
\cdot
R_n (\theta,\alpha)
\label{eq:RatioViaR}
\end{eqnarray}
with 
\begin{eqnarray}
R_n (\theta,\alpha)
&=&
\bE{ (1- \alpha \eta_{12}(\theta)) 
\left ( 1 + \widetilde Z (\theta,\alpha) \right )^{n-2} }
\nonumber \\
&=&
\bE{ (1-\eta_{12}(\theta) )
\left ( 1 + \widetilde Z (\theta,\alpha) \right )^{n-2} }
\label{eq:ExpressionForR} \\
& &
~  + (1-\alpha)
\bE{ \eta_{12}(\theta) \left ( 1 + \widetilde Z (\theta,\alpha) \right )^{n-2} } .
\nonumber
\end{eqnarray}

From the expression (\ref{eq:V}) it is plain that
\begin{equation} 
v(\theta;2K) \leq v(\theta; | K_1(\theta) \cup K_2(\theta)| ) \leq v(\theta;K) 
\label{eq:InequalitiesOnV}
\end{equation} 
with the lower (resp. upper) bound corresponding 
to $K_1(\theta) \cap K_2(\theta) = \emptyset$ 
(resp. $K_1(\theta) = K_2(\theta)$).

In the first term in (\ref{eq:ExpressionForR}),
the event $[ \eta_{12}(\theta) = 0 ]$ coincides with the event
$[ K_1(\theta) \cap K_2(\theta) = \emptyset ]$, 
in which case we have $| K_1(\theta) \cup K_2(\theta) | = 2K$ so that
\[
v(\theta; K_1(\theta) \cup K_2(\theta) ) - q(\theta)^2 
= v(\theta;2K) - q(\theta)^2 .
\]
We then conclude that
\begin{align}
& 
\bE{ (1-\eta_{12}(\theta)) \left ( 1 + \widetilde Z (\theta,\alpha) \right )^{n-2} } .
\nonumber \\
& = q(\theta) 
\left (
1 + \frac{ \alpha^2
\left (
v(\theta;2K) - q(\theta)^2
\right ) }
{ \left ( 1 - p(\theta,\alpha) \right )^2 }
\right )^{n-2}.
\label{eq:EqualityFirstTerm}
\end{align}

It is plain that
\begin{eqnarray}
\lefteqn{
 \left ( 1 - p(\theta,\alpha) \right )^2  - \alpha^2 q(\theta)^2 
 } & &
\nonumber \\
&=&
\left ( 1 - p(\theta,\alpha)  - \alpha q(\theta) \right )
\left ( 1 - p(\theta,\alpha)  + \alpha q(\theta) \right )
\nonumber \\
&=& 
\left ( 1 - \alpha \right )
\left ( 1 - \alpha + 2 \alpha q(\theta) \right ) > 0.
\end{eqnarray}
We also observe that
\[
v(\theta;2K) < q(\theta)^2
\]
by easy calculations based on the combinatorial expressions
for the quantities involved; details are left to the interested reader.
As a result, we have
\[
0 \leq 1 + \frac{ \alpha^2
\left (
v(\theta;2K) - q(\theta)^2
\right ) }
{ \left ( 1 - p(\theta,\alpha) \right )^2 }
\leq 1
\]
and the conclusion
\begin{eqnarray}
\bE{ (1-\eta_{12}(\theta)) \left ( 1 + \widetilde Z (\theta,\alpha) \right )^{n-2} } 
\leq q(\theta) 
\label{eq:BasicBound-}
\end{eqnarray}
follows.

As we turn to the second  term in (\ref{eq:ExpressionForR}),
it follows from (\ref{eq:InequalitiesOnV}) that
\begin{eqnarray}
\widetilde Z(\theta,\alpha)
&\leq&
\frac{ \alpha^2
\left (
v(\theta;K) - q(\theta)^2
\right ) }
{ \left ( 1 - p(\theta,\alpha) \right )^2 }
\nonumber \\
&=&
\frac{ \alpha^2
\left ( q(\theta) - q(\theta)^2 
\right ) }
{ \left ( 1 - p(\theta,\alpha) \right )^2 }
\nonumber \\
&=&
\frac{ \alpha^2 q(\theta) \left ( 1 - q(\theta) \right ) }
     { \left ( 1 - p(\theta,\alpha) \right )^2 }.
\end{eqnarray}

Using this deterministic bound we obtain
\begin{align}
& 
(1-\alpha) 
\bE{ \eta_{12}(\theta) \left ( 1 + \widetilde Z (\theta,\alpha) \right )^{n-2} } 
\nonumber \\
& \leq
(1-\alpha)
\bE{ \eta_{12}(\theta) }
\left (
1 + \frac{ \alpha^2 q(\theta) \left ( 1 - q(\theta) \right ) }
         { \left ( 1 - p(\theta,\alpha) \right )^2 }
\right)^{n-2}
\nonumber \\
& = 
(1-\alpha) (1-q(\theta))
\left (
1 + \frac{ \alpha q(\theta) p(\theta,\alpha) }
     { \left ( 1 - p(\theta,\alpha) \right )^2 }
\right)^{n-2}
\nonumber \\
& \leq 
(1-\alpha) (1-q(\theta)) \cdot R_n^\star(\theta,\alpha)
\label{eq:BasicBound+}
\end{align}
where we have set
\begin{equation}
R_n^\star(\theta,\alpha)
= 
e^{ (n-2) \frac{ \alpha q(\theta) p(\theta,\alpha) }
                { \left ( 1 - p(\theta,\alpha) \right )^2 }
}.
\label{eq:BasicBound2}
\end{equation}
Collecting (\ref{eq:BasicBound-}) and (\ref{eq:BasicBound+})
we obtain the key bound
\begin{eqnarray}
R_n (\theta,\alpha) \leq q(\theta) + (1-q(\theta)) R_n^\star(\theta,\alpha) .
\label{eq:BasicBound}
\end{eqnarray}
Later on we shall also have use for the quantity 
\begin{equation}
R^\circ_n (\theta,\alpha)
= e^{ \left ( 1 - p(\theta,\alpha) \right )^{-2} \cdot \alpha \log n } .
\label{eq:NewBoundOnExponent3}
\end{equation}

\section{A proof of Theorem \ref{thm:PartialVeryStrongZeroLaw+NodeIsolation}: The basic approach}
\label{sec:ProofPartialVeryStrongZeroLaw+NodeIsolation}

The proof of the zero-law of
Theorem \ref{thm:PartialVeryStrongZeroLaw+NodeIsolation} is developed
in the next three sections.
For the remainder of the paper, we consider fixed scalings $K,P: \mathbb{N}_0 \rightarrow \mathbb{N}_0$
and $\alpha: \mathbb{N}_0 \rightarrow [0,1]$ such that
(\ref{eq:ScalingCondition+VeryStrong}) holds
for some deviation function $\gamma : \mathbb{N}_0 \rightarrow \mathbb{R}$.
We also assume that
(\ref{eq:ConditionForVeryStrongZeroLaw+NodeIsolation}) holds.

From (\ref{eq:BoundViaSecondMoment}) 
the zero-law 
$\lim_{n \to \infty} \bP{I_n(\theta_n,\alpha_n) = 0} = 0$ 
will be established if we can show that
\begin{equation}
\liminf_{n \to \infty} 
\frac{ \left ( \bE{ I_n (\theta_n,\alpha_n)} \right )^2}
     { \bE{ I_n (\theta_1,\alpha_n) ^2} }
\geq 1.
\label{eq:ZeroLaw+NodeIsolation+Convergence2}
\end{equation}
In view of (\ref{eq:SecondMomentExpression}) this will be achieved if
the limiting statements
\begin{equation}
\lim_{n \to \infty} \bE{ I_n (\theta_n,\alpha_n) } = \infty
\label{eq:ZeroLaw+NodeIsolation+Convergence3}
\end{equation}
and
\begin{equation}
\limsup_{n \to \infty} \left( \frac{\bE{ \chi_{n,1} (\theta_n,\alpha_n)
\chi_{n,2} (\theta_n,\alpha_n) }}
     {\left (  \bE{ \chi_{n,1} (\theta_n,\alpha_n) } \right )^2 }
\right) \leq 1
\label{eq:ZeroLaw+NodeIsolation+Convergence4}
\end{equation}
both hold.

As the former holds by virtue of 
Lemma \ref{lem:VeryStrongZeroOneLaw+FirstMoment}
under (\ref{eq:ConditionForVeryStrongZeroLaw+NodeIsolation}),
it remains only to show the latter.
Using (\ref{eq:Limit=Zero}) we conclude from (\ref{eq:RatioViaR})
that establishing 
(\ref{eq:ZeroLaw+NodeIsolation+Convergence4}) is equivalent 
to showing
\begin{equation}
\limsup_{n\rightarrow \infty} 
R_n(\theta_n, \alpha_n) \leq 1 ,
\label{eq:ZeroLaw+NodeIsolation+Convergence5}
\end{equation}
and this will hold if we show the {\em stronger} inequality
\begin{equation}
\limsup_{n\rightarrow \infty} 
\left ( q(\theta_n) + (1-q(\theta_n)) R_n^\star(\theta_n,\alpha_n) \right )
\leq 1 .
\label{eq:StrongerVersionBasicBound}
\end{equation}

Under (\ref{eq:ScalingCondition+VeryStrong}), 
by the remarks made in the proof of Lemma 
\ref{lem:VeryStrongZeroOneLaw+FirstMoment},
we see that the exponent in $R_n^\star(\theta_n,\alpha_n)$ satisfies 
\begin{align}
&  (n-2) \frac{ \alpha_n q(\theta_n) \cdot p(\theta_n,\alpha_n) }
               { \left ( 1 - p(\theta_n,\alpha_n) \right )^2 }
\nonumber \\
& =
\frac{n-2}{n}
\frac{ \alpha_n q(\theta_n) \log n }
     { \left ( 1 - p(\theta_n,\alpha_n) \right )^2 }
-
\frac{n-2}{n}
\frac{ \alpha_n q(\theta_n) |\gamma_n| }
     { \left ( 1 - p(\theta_n,\alpha_n) \right )^2 }
\nonumber \\
& \leq
\left ( 1 - p(\theta_n,\alpha_n) \right )^{-2}
\cdot \alpha_n \log n
\label{eq:NewBoundOnExponent}
\end{align}
for $n=2,3, \ldots$ sufficiently large.
On that range this leads to the bound
\begin{equation}
R^\star_n(\theta_n,\alpha_n)
\leq R^\circ_n (\theta_n,\alpha_n)
\label{eq:NewBoundOnExponent2}
\end{equation}
as we recall (\ref{eq:NewBoundOnExponent3}).
Therefore, 
(\ref{eq:StrongerVersionBasicBound}) will hold 
if we show the stronger statement 
\begin{eqnarray}
\limsup_{n\rightarrow \infty}
\left ( 
q(\theta_n) + (1-q(\theta_n)) R^\circ_n (\theta_n,\alpha_n) \right )
\leq 1.
\label{eq:NewStrongerVersionBasicBound}
\end{eqnarray}

During the discussion we shall make use of the following two observations:
First the equality
\begin{eqnarray}
\lefteqn{
q(\theta_n) + (1-q(\theta_n))  R^\circ_n(\theta_n,\alpha_n)
} & &
\nonumber \\
&=&
1 + (1-q(\theta_n)) ( R^\circ_n(\theta_n,\alpha_n) - 1)
\label{eq:ImportantEquality}
\end{eqnarray}
holds for all $n=1,2, \ldots$.
Next, 
as already noted in the proof of Lemma \ref{lem:VeryStrongZeroOneLaw+FirstMoment},
condition
(\ref{eq:ConditionForVeryStrongZeroLaw+NodeIsolation}) yields 
$\gamma_n = -|\gamma_n|$ and  $|\gamma_{n} | \leq \log n$ eventually.
Thus, for $n=1,2, \ldots $ sufficiently large, {\em whenever} it happens that $\alpha_n > 0$,  we have
the bounds
\begin{eqnarray}
1 - q(\theta_n)
&=& 
\frac{1}{\alpha_n} \cdot \frac{ \log n - | \gamma_n | }{n}
\nonumber \\
&\leq& 
\frac{1}{\alpha_n} \cdot \frac{ \log n}{n}
%\label{eq:INEQUALITY=A} \\
\nonumber \\
&=& 
\frac{1}{\alpha_n \log n } \cdot \frac{ ( \log n )^2}{n}.
\label{eq:INEQUALITY=B}
\end{eqnarray}

\section{Along subsequences}
\label{AlongSubsequences}

Several cases need to be considered on the basis of the
behavior of the sequence
$\{ \alpha_n \log n, \ n=1,2, \ldots \}$ along subsequences.

\begin{lemma}
{\sl
Assume that along the subsequence $k \rightarrow n_k$,
the limit $\lim_{k \rightarrow \infty} \alpha_{n_k} \log n_k $ exists with
\begin{equation}
\lim_{k \rightarrow \infty} \alpha_{n_k} \log n_k = 0.
\label{eq:LimitAlog=0}
\end{equation}
Then, under (\ref{eq:ScalingCondition+VeryStrong}) 
with (\ref{eq:ConditionForVeryStrongZeroLaw+NodeIsolation})
we have both
\begin{equation}
\lim_{k \rightarrow \infty} R^\circ_{n_k} (\theta_{n_k},\alpha_{n_k} )= 1 ,
\label{eq:VeryStrongLimitR^+=1}
\end{equation}
and 
\begin{equation}
\lim_{k \rightarrow \infty} 
\left ( q(\theta_{n_k}) 
+ ( 1 - q(\theta_{n_k})  ) R^\circ_{n_k} (\theta_{n_k},\alpha_{n_k} ) 
\right ) = 1.
\label{eq:VeryStrongLimitDesiredR^+=1}
\end{equation}
}
\label{lem:VeryStrongLimit(1-q)R^+=1}
\end{lemma}

\myproof
Under the enforced assumptions, we have
\[
\lim_{k \rightarrow \infty}
\left ( 1 - p(\theta_{n_k},\alpha_{n_k}) \right )^{-2}
\cdot \alpha_{n_k} \log n_k 
= 0 
\]
as we recall (\ref{eq:Limit=Zero}).
The conclusion (\ref{eq:VeryStrongLimitR^+=1})   is then straightforward
from the expression (\ref{eq:NewBoundOnExponent3}), and 
(\ref{eq:VeryStrongLimitDesiredR^+=1}) follows upon using 
(\ref{eq:ImportantEquality}).
\myendpf

In this last step we had no information concerning
$\lim_{k \rightarrow \infty} q(\theta_{n_k}) $, hence the need for (\ref{eq:ImportantEquality})
in order to conclude (\ref{eq:VeryStrongLimitDesiredR^+=1}).

\begin{lemma}
{\sl
Assume that along the subsequence $k \rightarrow n_k$, the limit
$\lim_{k \rightarrow \infty} \alpha_{n_k} \log n_k$ exists in $(0,\infty)$.
Then, under (\ref{eq:ScalingCondition+VeryStrong}) 
with (\ref{eq:ConditionForVeryStrongZeroLaw+NodeIsolation})
we have both
\begin{equation}
\lim_{k \rightarrow \infty}
q(\theta_{n_k})  =1 
\label{eq:VeryStrong+Limit(1-q)=0}
\end{equation}
and
\begin{equation}
\lim_{k \rightarrow \infty}
\left ( 1 - q(\theta_{n_k}) \right )
R^\circ_{n_k} (\theta_{n_k},\alpha_{n_k} )
= 0,
\label{eq:VeryStrong+Limit(1-q)R^+=0}
\end{equation}
whence (\ref{eq:VeryStrongLimitDesiredR^+=1}) holds.
}
\label{lem:VeryStrongLimit(1-q)and(1-q)R^+=0}
\end{lemma}

\myproof
The condition $\lim_{k \rightarrow \infty} \alpha_{n_k} \log n_k > 0$
implies $\alpha_{n_k} > 0$ eventually. This together with condition
(\ref{eq:ConditionForVeryStrongZeroLaw+NodeIsolation}) allows us to use
(\ref{eq:INEQUALITY=B}) eventually along the subsequence $k \rightarrow n_k$.
Thus, for all $k=1,2, \ldots $ sufficiently large, we have
\begin{eqnarray}
1 - q(\theta_{n_k})
\leq
\frac{1}{ \alpha_{n_k} \log n_k }
\cdot 
\left (
\frac{ ( \log n_k )^2 }{ n_k }
\right ),
\label{eq:BOUND=1}
\end{eqnarray}
and the conclusion (\ref{eq:VeryStrong+Limit(1-q)=0}) immediately follows.
Finally, using (\ref{eq:Limit=Zero}) we get
\begin{equation}
\lim_{k \rightarrow \infty} R^\circ_{n_k} (\theta_{n_k},\alpha_{n_k} )
= 
e^{ \lim_{k \rightarrow \infty} \alpha_{n_k} \log n_k }
\label{eq:Intermediary1}
\end{equation}
where the limit is  finite by assumption, and the
conclusion (\ref{eq:VeryStrong+Limit(1-q)R^+=0}) follows
from  (\ref{eq:VeryStrong+Limit(1-q)=0}). The convergence
(\ref{eq:VeryStrongLimitDesiredR^+=1}) is now straightforward.
\myendpf

\begin{lemma}
{\sl
Assume that along the subsequence $k \rightarrow n_k$,
the limit $\lim_{k \rightarrow \infty} \alpha_{n_k} \log n_k $ exists with
\begin{equation}
\lim_{k \rightarrow \infty} \alpha_{n_k} \log n_k = \infty.
\label{eq:LimitAlog=Infinity}
\end{equation}
Then, under (\ref{eq:ScalingCondition+VeryStrong}) 
with (\ref{eq:ConditionForVeryStrongZeroLaw+NodeIsolation})
we still have (\ref{eq:VeryStrong+Limit(1-q)=0}) whereas both
(\ref{eq:VeryStrongLimitDesiredR^+=1}) and (\ref{eq:VeryStrong+Limit(1-q)R^+=0})
hold provided the additional condition
\begin{equation}
\limsup_{k \rightarrow \infty} \alpha_{n_k}  < 1 
\label{eq:AdditionalCondition}
\end{equation}
is enforced.
}
\label{lem:AgainVeryStrongLimit(1-q)R^+=1}
\end{lemma}

\myproof
It is plain that (\ref{eq:VeryStrong+Limit(1-q)=0}) still holds under the condition
$\lim_{k \rightarrow \infty}  \alpha_{n_k} \log n_k = \infty$ since the bound (\ref{eq:BOUND=1}) 
is valid here as well since $\alpha_{n_k}  > 0$ eventually for all $k=1,2, \ldots$ sufficiently large.
In order to justify (\ref{eq:VeryStrong+Limit(1-q)R^+=0}) under
the additional condition (\ref{eq:AdditionalCondition}) we argue as follows:
Consider $k =1,2, \ldots $ sufficiently large so that (\ref{eq:BOUND=1}) holds.
We have
\begin{eqnarray}
\lefteqn{
(1-q(\theta_{n_k}) ) R^\circ_{n_k} ( \theta_{n_k}, \alpha_{n_k} )
} & &
\nonumber \\
&\leq&
\frac{1}{ \alpha_{n_k} \log n_k }
\cdot 
\left (
\frac{ ( \log n_k )^2 }{ n_k }
\right )
R^\circ_{n_k} ( \theta_{n_k}, \alpha_{n_k})
\nonumber \\
&=&
\frac{1}{ \alpha_{n_k} \log n_k }
\cdot ( \log n_k )^2 
\cdot
\frac{ R^\circ_{n_k} ( \theta_{n_k}, \alpha_{n_k}) }{ n_k}
\label{eq:IntermediaryRelation=A}
\end{eqnarray}
with
\begin{eqnarray}
\lefteqn{
\frac{ R^\circ_{n_k} ( \theta_{n_k}, \alpha_{n_k}) }{ n_k}
}
& & 
\nonumber \\
&=& \frac{1}{n_k} \cdot
e^{ (1-p(\theta_{n_k},\alpha_{n_k} ))^{-2} \cdot \alpha_{n_k} \log n_k }
\nonumber \\
&=&
e^{ \left ( -1 + (1-p(\theta_{n_k},\alpha_{n_k} ))^{-2} \cdot \alpha_{n_k}  \right ) \log n_k }.
\nonumber 
\end{eqnarray}

By virtue of (\ref{eq:Limit=Zero}), we find that
\begin{eqnarray}
\lefteqn{
\limsup_{k \rightarrow \infty} 
\left ( -1 + (1-p(\theta_{n_k},\alpha_{n_k} ))^{-2} \cdot \alpha_{n_k}  \right )
} & &
\nonumber \\
&=&
-1 + \limsup_{k \rightarrow \infty} 
 \frac{ \alpha_{n_k} }{(1-p(\theta_{n_k},\alpha_{n_k} ))^2 }
\nonumber \\
&=&
-1 + \left ( \limsup_{k \rightarrow \infty}  \alpha_{n_k}  \right )  < 0,
\end{eqnarray}
under the additional condition (\ref{eq:AdditionalCondition}), whence
\begin{equation}
\lim_{k \rightarrow \infty} 
\left ( ( \log n_k )^2 
\cdot
\frac{ R^\circ_{n_k} ( \theta_{n_k}, \alpha_{n_k}) }{ n_k} 
\right )= 0.
\label{eq:LimitFunnyTerm}
\end{equation}

Let $k$ to infinity in (\ref{eq:IntermediaryRelation=A}): 
The validity of (\ref{eq:VeryStrong+Limit(1-q)R^+=0}) now follows by appealing to
(\ref{eq:LimitAlog=Infinity}) and (\ref{eq:LimitFunnyTerm}).
Here as well the convergence
(\ref{eq:VeryStrongLimitDesiredR^+=1}) is straightforward.
\myendpf

In summary, under their specific assumptions, each of 
Lemma \ref{lem:VeryStrongLimit(1-q)R^+=1},
Lemma \ref{lem:VeryStrongLimit(1-q)and(1-q)R^+=0} and
Lemma \ref{lem:AgainVeryStrongLimit(1-q)R^+=1}
ensures that (\ref{eq:VeryStrongLimitDesiredR^+=1}) holds, hence
$\limsup_{k \rightarrow \infty} R_{n_k} ( \theta_{n_k},\alpha_{n_k} ) \leq 1$
and the conclusion
\[
\lim_{k \rightarrow \infty} P_{n_k} ( \theta_{n_k},\alpha_{n_k} ) = 0
\]
follows.

\section{Completing the proof of Theorem \ref{thm:PartialVeryStrongZeroLaw+NodeIsolation}}
\label{sec:CompleteProof}

The proof of Theorem \ref{thm:PartialVeryStrongZeroLaw+NodeIsolation} relies on the Subsequence Principle:
For any {\em arbitrary} subsequence $k \rightarrow n_k$, we shall show that there exists 
a further subsequence $\ell \rightarrow k_\ell$ such that
\begin{equation}
\lim_{\ell \rightarrow \infty} 
P_{n_{k_\ell}} ( \theta_{n_{k_\ell}} ,   \alpha_{n_{k_\ell}} ) = 0.
\label{eq:SubsequenceLimit}
\end{equation}
It is well known that this implies $\lim_{n \rightarrow \infty} P_n(\theta_n, \alpha_n) = 0$.

If 
\[
\limsup_{n \rightarrow \infty} \alpha_n \log n < \infty,
\]
then $\limsup_{k \rightarrow \infty} \alpha_{n_k} \log {n_k} < \infty$ as well,
and there exists a subsequence $\ell \rightarrow k_\ell$ such that
\[
\lim_{\ell \rightarrow \infty}  \alpha_{n_{k_\ell}} \log n_{k_\ell}
=
\limsup_{k \rightarrow \infty} \alpha_{n_k} \log {n_k} 
\]
When $ \limsup_{k \rightarrow \infty} \alpha_{n_k} \log {n_k}  = 0$, we invoke 
Lemma  \ref{lem:VeryStrongLimit(1-q)R^+=1} (applied to the subsequence $\ell \rightarrow n_{k_\ell}$)
to conclude that (\ref{eq:SubsequenceLimit}) holds.
On the other hand, if $ \limsup_{k \rightarrow \infty} \alpha_{n_k} \log {n_k}  $ is an element of $(0,\infty)$ we also
conclude to (\ref{eq:SubsequenceLimit}) by appealing to 
Lemma \ref{lem:VeryStrongLimit(1-q)and(1-q)R^+=0} (applied to the subsequence $\ell \rightarrow n_{k_\ell}$).

If 
\[
\limsup_{n \rightarrow \infty} \alpha_n \log n = \infty,
\]
then there are two possibilities:

(i) If $\limsup_{k \rightarrow \infty} \alpha_{n_k} \log {n_k} <  \infty$,
then the earlier analysis applies unchanged and leads to the existence
of a subsequence $\ell \rightarrow k_\ell$ such that  (\ref{eq:SubsequenceLimit}) holds.

(ii) If $\limsup_{k \rightarrow \infty} \alpha_{n_k} \log {n_k} = \infty$,
 there exists at least one subsequence $\ell \rightarrow k_\ell$ such that
\[
\lim_{\ell \rightarrow \infty}  \alpha_{n_{k_\ell}} \log n_{k_\ell} = \infty
\]
On the other hand, the condition $\limsup_{n \rightarrow \infty} \alpha_n < 1$ implies
$\limsup_{\ell \rightarrow \infty}  \alpha_{n_{k_\ell}} < 1 $, and 
Lemma \ref{lem:AgainVeryStrongLimit(1-q)R^+=1} (applied to the subsequence $\ell \rightarrow n_{k_\ell}$)
ensures that (\ref{eq:SubsequenceLimit}) holds.
 
\section*{Acknowledgment}
This work was supported by NSF Grant CCF-1217997.
The paper was completed during the academic year 2014-2015  while A.M. Makowski 
was a Visiting Professor with the Department of Statistics of the Hebrew University of Jerusalem 
with the support of a fellowship from the Lady Davis Trust.

\bibliographystyle{IEEE}

%\begin{thebibliography}{10}
%\bibliography{../main}

\section{Appendix: Proof of (\ref{eq:FirstMomentCHI}) and (\ref{eq:CrossMoment})}
\label{sec:Appendix}

Fix $n=2,3, \ldots$, positive integers $K$ and $P$ such that $K < P$,
and $\alpha$ in $[0,1]$.
From (\ref{eq:Defn+CHI}) we note that
\[
\chi_{n,j}(\theta,\alpha)
= 
\left (
1 - B_{12}(\alpha) \eta_{12}(\theta)
\right )
\cdot
\prod_{\ell=3}^n 
\left (
1 - B_{j\ell}(\alpha) \eta_{j\ell}(\theta)
\right )
\]
for $j=1,2$.

It follows from 
(\ref{eq:Probab_Keyring_Does_Not_Intersect_R}) that
the indicator rvs $\eta_{12}(\theta), \ldots , \eta_{1n}(\theta)$ are i.i.d.
$\{0,1\}$ rvs, so that the rvs $\eta_{12}(\theta), \ldots , \eta_{1n}(\theta) , B_{12}(\alpha), \ldots , B_{1n}(\alpha)$ are 
mutually independent under the enforced independence assumptions.
Therefore, 
\begin{eqnarray}
\bE{ \chi_{n,1}(\theta,\alpha) }
&=&
\bE{ 
\prod_{\ell=2}^n 
\left (
1 - B_{1\ell}(\alpha) \eta_{1\ell}(\theta)
\right )
}
\nonumber \\
&=&
\prod_{\ell=2}^n 
\bE{ 
1 - B_{1\ell}(\alpha) \eta_{1\ell}(\theta)
}
\nonumber \\
&=&
\prod_{\ell=2}^n 
\left (
1 - \bE{ B_{1\ell}(\alpha)} \bE{ \eta_{1\ell}(\theta) }
\right )
\nonumber \\
&=&
\prod_{\ell=2}^n 
\left ( 1 - \alpha (1- q(\theta) )\right ),
\end{eqnarray}
and we obtain the expression (\ref{eq:FirstMomentCHI}).

Next, we observe that
\[
 \chi_{n,1} (\theta,\alpha) \chi_{n,2} (\theta,\alpha) 
=  \left ( 1 - B_{12}(\alpha) \eta_{12}(\theta) \right )
\cdot
\left ( \ldots \right )
\]
with
\[
\ldots 
=
\prod_{\ell=3}^n 
\left (
1 - B_{1\ell}(\alpha) \eta_{1\ell}(\theta)
\right )
\left (
1 - B_{2\ell}(\alpha) \eta_{2\ell}(\theta)
\right ).
\]
Upon conditioning with respect to the rvs 
$K_1(\theta), \ldots , K_n(\theta)$, we get
\begin{align}
& 
\bE{ \chi_{n,1} (\theta,\alpha) \chi_{n,2} (\theta,\alpha) }
\nonumber \\
&=
\bE{ 
\left ( 1 - \alpha \eta_{12}(\theta) \right )
\cdot
\prod_{\ell=3}^n 
\left ( 1 - \alpha \eta_{1\ell}(\theta) \right )
\left ( 1 - \alpha \eta_{2\ell}(\theta) \right )
}
\nonumber
\end{align}
under the enforced independence assumptions.
It is also easy to check that {\em conditionally} on $K_1(\theta)$ and $K_2(\theta)$,
the $n-2$ pairs of rvs $(\eta_{13}(\theta), \eta_{23}(\theta) ), \ldots , (\eta_{1n}(\theta), \eta_{2n}(\theta) )$ are
mutually independent, whence
\begin{align}
& 
\bE{
\prod_{\ell=3}^n 
\left ( 1 - \alpha \eta_{1\ell}(\theta) \right )
\left ( 1 - \alpha \eta_{2\ell}(\theta) \right )
\Bigg | K_1(\theta), K_2(\theta) }
\nonumber \\
&=
\prod_{\ell=3}^n 
\bE{
\left ( 1 - \alpha \eta_{1\ell}(\theta) \right )
\left ( 1 - \alpha \eta_{2\ell}(\theta) \right )
| K_1(\theta), K_2(\theta) } .
\nonumber
\end{align}

For each $\ell =3, \ldots , n$ and $j=1,2$, it is plain that
\begin{eqnarray}
\lefteqn{ \bE{ \eta_{j\ell}(\theta) | K_1(\theta), K_2(\theta) } }  & &
\nonumber \\
&=& 
\bE{ \1{ K_j(\theta) \cap K_\ell(\theta) \neq \emptyset  } | K_i(\theta) }
\nonumber \\
&=& 1 - q(\theta).
\end{eqnarray}
Therefore,
\begin{align}
& 
\bE{ 
\left ( 1 - \alpha \eta_{1\ell}(\theta) \right )
\left ( 1 - \alpha \eta_{2\ell}(\theta) \right ) 
| K_1(\theta), K_2(\theta) }
\nonumber \\
&=
1 - 2\alpha (1-q(\theta) ) + \alpha^2  \bE{ \eta_{1\ell}(\theta) \eta_{2\ell}(\theta)  | K_1(\theta), K_2(\theta) }
\nonumber \\
&=
\left ( 1 - \alpha (1-q(\theta) ) \right )^2
\nonumber \\
&~~ 
+ \alpha^2
\left ( 
\bE{ \eta_{1\ell}(\theta) \eta_{2\ell}(\theta)  | K_1(\theta), K_2(\theta) }
-
(1 -q(\theta) )^2
\right )
\nonumber
\end{align}
by a completion-of-square argument.
With the product rv $\eta_{1\ell}(\theta) \eta_{2\ell}(\theta) $ given by
\begin{eqnarray}
\lefteqn{
\left ( 1 - \1{ K_1(\theta) \cap K_\ell (\theta) = \emptyset } \right ) 
\left ( 1 - \1{ K_2(\theta) \cap K_\ell (\theta) = \emptyset } \right ) 
} & &
\nonumber \\
&=&
1 - \1{ K_1(\theta) \cap K_\ell (\theta) = \emptyset }  - \1{ K_2(\theta) \cap K_\ell (\theta) = \emptyset } 
\nonumber \\
& & ~+ \1{ K_1(\theta) \cap K_\ell (\theta) = \emptyset } \1{ K_2(\theta) \cap K_\ell (\theta) = \emptyset } 
\nonumber \\
&=&
1 - \1{ K_1(\theta) \cap K_\ell (\theta) = \emptyset }  - \1{ K_2(\theta) \cap K_\ell (\theta) = \emptyset } 
\nonumber \\
& & ~+ \1{ ( K_1(\theta) \cup K_2(\theta) )  \cap K_\ell (\theta) = \emptyset } ,
\nonumber
\end{eqnarray}
it follows that
\begin{eqnarray}
\lefteqn{
 \bE{ \eta_{1\ell}(\theta) \eta_{2\ell}(\theta)  | K_1(\theta), K_2(\theta) }
 } & &
 \nonumber \\
 &=&
 1 - 2 q(\theta)
 + v(\theta; K_1(\theta) \cup K_2(\theta)  )
\end{eqnarray}
so that
\begin{eqnarray}
\lefteqn{
 \bE{ \eta_{1\ell}(\theta) \eta_{2\ell}(\theta)  | K_1(\theta), K_2(\theta) } 
 - (1 -q(\theta) )^2
 } & &
 \nonumber \\
 &=&
 1 - 2 q(\theta)
 + v(\theta; K_1(\theta) \cup K_2(\theta)  )
  - (1 -q(\theta) )^2
  \nonumber \\
  &=& 
 v(\theta; K_1(\theta) \cup K_2(\theta)  ) - q(\theta)^2.
\end{eqnarray}
Collecting terms we conclude that
\begin{align}
& 
\bE{ 
\left ( 1 - \alpha \eta_{1\ell}(\theta) \right )
\left ( 1 - \alpha \eta_{2\ell}(\theta) \right ) 
| K_1(\theta), K_2(\theta) }
\nonumber \\
&=
\left ( 1 - p(\theta,\alpha) \right )^2
+ \alpha^2
\left (
 v(\theta; K_1(\theta) \cup K_2(\theta)  ) - q(\theta)^2
\right ).
\nonumber
\end{align}
Upon substitution into earlier expressions, we now obtain
\begin{align}
& 
\bE{
\prod_{\ell=3}^n 
\left ( 1 - \alpha \eta_{1\ell}(\theta) \right )
\left ( 1 - \alpha \eta_{2\ell}(\theta) \right )
\Bigg | K_1(\theta), K_2(\theta) }
\nonumber \\
&=
Z(\theta;\alpha)^{n-2}
\nonumber
\end{align}
with the rv $Z(\theta;\alpha)$ is given by
\begin{eqnarray}
\lefteqn{
Z(\theta;\alpha)
} & &
\label{eq:ZZZ}  \\
&=&
\left ( 1 - p(\theta,\alpha) \right )^2
+ \alpha^2
\left (
 v(\theta; K_1(\theta) \cup K_2(\theta)  ) - q(\theta)^2
\right )
\nonumber
\end{eqnarray}
Finally,
\[
\bE{ \chi_{n,1} (\theta,\alpha) \chi_{n,2} (\theta,\alpha) }
=
\bE{ \left ( 1 - \alpha \eta_{12}(\theta) \right ) Z(\theta;\alpha)^{n-2} }.
\]
and we note that the rv $Z(\theta;\alpha)$ given at (\ref{eq:ZZZ}) coincides with the rv
given through the expressions (\ref{eq:Z})-(\ref{eq:TildeZ}).
\myendpf

\end{document}